# A Leray-Trudinger Inequality

G. Psaradakis[a,*], D. Spector[a,b]

[a]*Department of Mathematics, Technion - Israel Institute of Technology, Haifa 32000, Israel*
[b]*Department of Applied Mathematics, National Chiao Tung University, Hsinchu 30010, Taiwan*

## Abstract

We consider a multidimensional version of an inequality due to Leray as a substitute for Hardy's inequality in the case $p = n \geq 2$. In this paper we provide an optimal Sobolev-type improvement of this substitute, analogous to the corresponding improvements obtained for $p = 2 < n$ in S. Filippas, A. Tertikas, Optimizing improved Hardy inequalities, J. Funct. Anal. 192 (1) (2002) 186–233, and for $p > n \geq 1$ in G. Psaradakis, An optimal Hardy-Morrey inequality, Calc. Var. Partial Differential Equations 45 (3-4) (2012) 421–441.

**MSC 2010**: 46E35; 26D15; 35J92

**Key Words**: Hardy inequality, Leray potential, borderline Sobolev embedding

## 1. Introduction

Let $\Omega \subset \mathbb{R}^n$ be open and connected. A multidimensional version of Hardy's inequality asserts that for $n > 2$ one has

$$\int_\Omega |\nabla u|^2 \mathrm{d}x \geq \left(\frac{n-2}{2}\right)^2 \int_\Omega \frac{|u|^2}{|x|^2}\mathrm{d}x \quad \text{for all } u \in C_c^\infty(\Omega), \tag{1.1}$$

with the best possible constant in case $0 \in \Omega$ (see [28] and [22]). The precise value of the constant in (1.1) plays, for example, a crucial role in the analysis of solutions of heat equations with potentials having critical point singularities (see [10], [11] and [32]).

If we define the Hardy difference

$$I[u;\Omega] := \int_\Omega |\nabla u|^2 \mathrm{d}x - \left(\frac{n-2}{2}\right)^2 \int_\Omega \frac{|u|^2}{|x|^2}\mathrm{d}x,$$

then (1.1) reads $I[u;\Omega] \geq 0$ for all $u \in C_c^\infty(\Omega)$. Actually, it is well known that $I[u;\Omega] > 0$ if $u \in W_0^{1,2}(\Omega) \setminus \{0\}$, which suggests the possibility of improving the inequality (1.1) in the form of lower bounds for $I[u;\Omega]$. While for $\Omega = \mathbb{R}^n$ it has been shown that additional correction terms cannot be added (see for example [16], [19] & [14]), if $\Omega$ has finite volume such an improvement is possible. The following subcritical Sobolev improvement of the Hardy inequality (1.1) is due to Brezis and

*Corresponding author

*Email addresses:* georgios@tx.technion.ac.il (G. Psaradakis),
dspector@tx.technion.ac.il (D. Spector)

Vazquez [10]: *If $\Omega$ has finite volume, then for any $1 \leq q < 2^*$, there exists $C_{n,q} > 0$ depending only on $n$ and $q$, such that*

$$\left(I[u;\Omega]\right)^{1/2} \geq \frac{C_{n,q}}{\mathrm{vol}(\Omega)^{1/q-1/2^*}} \left(\int_\Omega |u|^q dx\right)^{1/q} \quad \text{for all } u \in C_c^\infty(\Omega). \tag{1.2}$$

Here, $2^* := 2n/(n-2)$ is the Sobolev critical exponent. We recall that for any domain $\Omega$ with finite volume, $2^*$ is the largest value of $q$ for which one has the existence of a constant $S_n > 0$ depending only on $n$ such that

$$\left(\int_\Omega |\nabla u|^2 dx\right)^{1/2} \geq \frac{S_n}{\mathrm{vol}(\Omega)^{1/q-1/2^*}} \left(\int_\Omega |u|^q dx\right)^{1/q} \quad \text{for all } u \in C_c^\infty(\Omega). \tag{1.3}$$

It was a question in [10] whether there is a further improvement of inequality (1.2). Filippas and Tertikas showed in [16] that though (1.2) fails for $q = 2^*$, introducing a logarithmic relaxation one can have a critical Sobolev improvement to Hardy's inequality. Their result is as follows.

HARDY-SOBOLEV INEQUALITY: *Let $\Omega$ be a bounded domain in $\mathbb{R}^n$; $n \geq 3$, containing the origin. Then there exists a constant $C_n > 0$ depending only on $n$, such that*

$$\left(I[u,\Omega]\right)^{1/2} \geq C_n \left(\int_\Omega \left(|u|X^{1-1/n}(|x|/R_\Omega)\right)^{2^*} dx\right)^{1/2^*} \quad \text{for all } u \in C_c^\infty(\Omega), \tag{1.4}$$

*where $R_\Omega := \sup_{x \in \Omega} |x|$ and $X(t) := (1-\log t)^{-1}$; $t \in (0,1]$. Moreover, the exponent $1 - 1/n$ on $X$ is optimal in the sense that it cannot be decreased.*

A non-trivial substitute of (1.1) in the case $n = 2$ is due to Leray [24], who used it in the study of two dimensional viscous flows. More generally, in analogy with versions of Hardy's inequality for $p \neq n$, it has been extended to $p = n \geq 2$ by ([3], [5] & [7]), and can be stated as follows: *If $\Omega$ is a bounded domain in $\mathbb{R}^n$; $n \geq 2$, then*

$$\int_\Omega |\nabla u|^n dx \geq \left(\frac{n-1}{n}\right)^n \int_\Omega \frac{|u|^n}{|x|^n} X^n(|x|/R_\Omega) dx \quad \text{for all } u \in C_c^\infty(\Omega \setminus \{0\}), \tag{1.5}$$

*with the best possible constant in case $0 \in \Omega$.*

If we define the Leray difference

$$I_n[u;\Omega] := \int_\Omega |\nabla u|^n dx - \left(\frac{n-1}{n}\right)^n \int_\Omega \frac{|u|^n}{|x|^n} X^n(|x|/R_\Omega) dx,$$

we have again that $I_n[u;\Omega] > 0$ for $u \in W_0^{1,n}(\Omega) \setminus \{0\}$ (see [8]). In analogy with the results of [10] and [16], it is then natural to ask whether one can make subcritical or critical Sobolev improvements to the inequality (1.5) via finding lower bounds for $I_n[u;\Omega]$. Here we enter into another type of criticality, where the Sobolev critical exponent is formally $+\infty$. The full understanding of what should be the analog of (1.3) when $n = 2$, and more generally when $p = n \geq 2$ was given by Trudinger in [31] (see also Peetre [27]), who proved the following result: *If $\Omega$ is a domain in* $\mathbb{R}^n$; $n \geq 2$, *having finite volume* $|\Omega|$, *then there exist positive constants* $a_n$ *and* $b_n$, *depending only on n, such that*

$$\int_\Omega e^{a_n |u|^{n/(n-1)}} dx \leq b_n \text{vol}(\Omega) \quad \text{for all } u \in W_0^{1,n}(\Omega) \text{ satisfying } \int_\Omega |\nabla u|^n dx \leq 1, \tag{1.6}$$

*with the optimal exponent on* $|u|$ *(cannot be increased).* For further information, historical notes, various extensions, sharp constants and applications of this important inequality we refer to [27], [2]-§3.8, [25], [6]-§2.15 & §2.16 and references therein.

Therefore, in view of the Hardy-Sobolev inequality (1.4), one wonders whether we have some exponential integrability of functions satisfying $I_n[u;\Omega] \leq 1$. In analogy with the results in [16], we show that though the direct combination fails to hold, one can obtain a critical Leray-Trudinger inequality with the introduction of a logarithmic correction. Our result is as follows.

**Theorem 1.1.** LERAY-TRUDINGER INEQUALITY: *Let* $\Omega$ *be a bounded domain in* $\mathbb{R}^n$; $n \geq 2$, *containing the origin. For any* $\varepsilon > 0$ *there exist positive constants* $A_{n,\varepsilon}$ *depending only on* $n, \varepsilon$, *and* $B_n$ *depending only on n, such that*

$$\int_\Omega e^{A_{n,\varepsilon} [|u(x)| X^\varepsilon (|x|/R_\Omega)]^{n/(n-1)}} dx \leq B_n \text{vol}(\Omega) \quad \text{for all } u \in C_c^\infty(\Omega \setminus \{0\}) \text{ satisfying } I_n[u;\Omega] \leq 1. \tag{1.7}$$

*Moreover, such an estimate fails for* $\varepsilon = 0$.

An interesting point to note here is that the exponent of the logarithmic correction can be chosen in the open interval $(0, +\infty)$, which is in stark contrast to (1.4) and also the case $p > n$ (see [26]), where the exponent lies in a closed interval. Note also that (1.2) for bounded domains can be obtained from (1.4) through a simple use of Hölder's inequality. Similar arguments enable one to deduce analogous subcritical results in our setting.

A related two dimensional result is in [30], where the author improves Moser's inequality (Trudinger's inequality (1.6) having optimal $a_2$ constant). This is another interesting perspective, where one begins with a critical Sobolev inequality with best constant and asks whether one can make "subcritical Hardy improvements". Our result focuses instead on the optimal allowed singularity of the potential, and moreover, it is valid in any dimension $n \geq 2$. Another related result is in [33], where the Hardy inequality involving distance to the boundary of a disc in $\mathbb{R}^2$ is considered instead of (1.5). Finally, let us mention [4], where an improvement of Moser's inequality for $n = 2$ was proven.

For the proof of (1.7) we follow closely Trudinger's original proof (see also [18]-§7.8) taking into account the corresponding ground state transform. As in [10], the ground state is a solution to the Euler-Lagrange equation associated to critical points of the best constant problem related to inequality (1.5). It is by now well understood that the exponential integrability of functions in $W_0^{1,n}(\Omega)$ rests on the following $L^q$ estimates: There exists a constant $C > 0$ depending only on $n$ and $\Omega$, such that

$$\|u\|_{L^q(\Omega)} \leq C q^{1/q + 1 - 1/n} \|\nabla u\|_{L^n(\Omega)},$$

for all $q$ sufficiently large. Our $L^q$ estimates read as follows (see Proposition 3.1): There exists a constant $C > 0$ depending only on $n$ and $\Omega$, such that for any $\varepsilon > 0$ and all $u \in C_c^\infty(\Omega \setminus \{0\})$

$$\|uX^\varepsilon\|_{L^q(\Omega)} \leq C \frac{q^{1/q + 1 - 1/n}}{\varepsilon} \left( I_n[u;\Omega] \right)^{1/n},$$

for all $q$ sufficiently large. To prove that (1.7) does not hold when $\varepsilon = 0$, the following optimal homogeneous improvement to (1.5) (found in [8]) plays a significant role: *In a bounded domain $\Omega$ of $\mathbb{R}^n$ containing the origin we have*

$$I_n[u;\Omega] \geq \frac{1}{2}\left(\frac{n-1}{n}\right)^{n-1} \int_\Omega \frac{|u|^n}{|x|^n} X^n(|x|/R_\Omega) X^2(X(|x|/R_\Omega)) \mathrm{d}x \quad \text{for all } u \in C_c^\infty(\Omega \setminus \{0\}), \quad (1.8)$$

*with the best possible constant. In addition, the exponent 2 on $X(X)$ cannot be decreased.* Assuming that (1.7) is true with $\varepsilon = 0$, we are able to show that we can improve the exponent 2 on $X(X)$, a contradiction.

For other directions in strengthening the inequality (1.5), we refer to [12]. In strengthening Trudinger's inequality (1.6), we refer to [13]. For the combination of Hardy's inequality and the Sobolev or Morrey inequality in the case $p > n$, see [26]. For Hardy-Sobolev inequalities with the weight being the distance to the boundary and $2 \leq p < n$ we refer to [9], [15] and [17].

## 2. A Hardy type inequality when $p = n$

This section is a discussion on the extension of Leray's inequality in any dimension $n \geq 2$. This plays the role of Hardy's inequality in the case $p = n \geq 2$.

In what follows, $B_R$ stands for an open ball in $\mathbb{R}^n$ having radius $R > 0$ and center at 0. The volume of $B_1$ is denoted by $\omega_n$. Also, by $\Omega$ we denote a bounded domain (open, connected set) in $\mathbb{R}^n$; $n \geq 2$. We set

$$R_\Omega := \sup_{x \in \Omega} |x|.$$

We note that if $\Omega$ contains the origin then trivially $\Omega \subseteq B_{R_\Omega}$, so that $\mathrm{vol}(\Omega) \leq \mathrm{vol}(B_{R_\Omega}) = \omega_n R_\Omega^n$. We also define the auxiliary function

$$X(t) := (1 - \log t)^{-1}, \quad \text{whenever } t \in (0,1].$$

This function is strictly increasing with $X(0+) = 0$ and $X(1) = 1$. Moreover, the following differentiation rule can be easily checked whenever $\gamma \in \mathbb{R}$

$$[X^\gamma(t)]' = \frac{\gamma}{t} X^{\gamma+1}(t); \quad t \in (0,1].$$

The following lemma is a weighted version of inequality (1.5). It follows by the choice $p = k = n$ in Lemma 3.2 of [7], though we give another proof for the convenience of the reader.

**Lemma 2.1.** *For all $\alpha \neq 1$ and any $u \in C_c^\infty(\Omega \setminus \{0\})$ we have*

$$\int_\Omega |\nabla u|^n X^{\alpha - n}(|x|/R_\Omega) \mathrm{d}x \geq \left|\frac{\alpha - 1}{n}\right|^n \int_\Omega \frac{|u|^n}{|x|^n} X^\alpha(|x|/R_\Omega) \mathrm{d}x.$$

**Proof.** It suffices to prove it in the case $R_\Omega = 1$. The general case then follows by a change of variables and density arguments. Integrating by parts we get

$$-\int_\Omega \nabla |u|^n \cdot \{|x|^{-n} X^{\alpha-1}(|x|)x\} dx = \int_\Omega |u|^n \mathrm{div}\{|x|^{-n} X^{\alpha-1}(|x|)x\} dx$$
$$= (\alpha - 1)\int_\Omega |u|^n |x|^{-n} X^\alpha(|x|) dx.$$

Thus we conclude

$$\int_\Omega \frac{|u|^{n-1}|\nabla u|}{|x|^{n-1}} X^{\alpha-1}(|x|) dx \geq \frac{|\alpha-1|}{n} \int_\Omega \frac{|u|^n}{|x|^n} X^\alpha(|x|) dx. \tag{2.1}$$

The left hand side can be written as follows

$$\int_\Omega \frac{|u|^{n-1}|\nabla u|}{|x|^{n-1}} X^{\alpha-1}(|x|) dx = \int_\Omega \left\{|\nabla u| X^{\alpha/n-1}(|x|)\right\}\left\{\frac{|u|^{n-1}}{|x|^{n-1}} X^{\alpha(n-1)/n}(|x|)\right\} dx$$
$$\leq \left(\int_\Omega |\nabla u|^n X^{\alpha-n}(|x|) dx\right)^{1/n} \left(\int_\Omega \frac{|u|^n}{|x|^n} X^\alpha(|x|) dx\right)^{1-1/n},$$

by Hölder's inequality. Combining this with (2.1), rearranging and taking the $n$th power of both sides, the result is demonstrated. ∎

**Remark 2.2.** The classic multidimensional Hardy inequality

$$\int_{\mathbb{R}^n} |\nabla u|^p dx \geq \left|\frac{p-n}{n}\right|^n \int_{\mathbb{R}^n} \frac{|u|^p}{|x|^p} dx; \quad 1 \leq p \neq n, \tag{2.2}$$

is valid for all $u \in C_c^\infty(\mathbb{R}^n \setminus \{0\})$. The constant is well known to be the optimal one. For $\alpha = n$ in Lemma 2.1 we obtain (see also [3], [5] and [29, Lemma 17.4])

$$\int_\Omega |\nabla u|^n dx \geq \left(\frac{n-1}{n}\right)^n \int_\Omega \frac{|u|^n}{|x|^n} X^n(|x|/R_\Omega) dx, \tag{2.3}$$

for all $u \in C_c^\infty(\Omega \setminus \{0\})$. If $0 \in \Omega$, then the constant is known to be the optimal one (see [5] and [7]). Thus, (2.3) may be considered as a "substitute" of the Hardy inequality (2.2) in case $p = n$, which is valid in bounded domains containing the origin. For other substitutes (even in $\mathbb{R}^n$) see [20].

**Definition 2.3.** Whenever $u \in C_c^\infty(\Omega \setminus \{0\})$ we set

$$I_n[u;\Omega] := \int_\Omega |\nabla u|^n dx - \left(\frac{n-1}{n}\right)^n \int_\Omega \frac{|u|^n}{|x|^n} X^n(|x|/R_\Omega) dx \ (\geq 0).$$

**Remark 2.4.** The function $u_0 : \Omega \setminus \{0\} \mapsto (1,\infty)$ defined by

$$u_0(x) := X^{-1+1/n}(|x|/R_\Omega),$$

is such that

$$-\mathrm{div}(|\nabla u_0|^{n-2}\nabla u_0) - \left(\frac{n-1}{n}\right)^n \frac{X^n(|x|/R_\Omega)}{|x|^n} u_0^{n-1} = 0 \quad \text{in } \Omega \setminus \{0\}.$$

**Definition 2.5.** Whenever $v \in C_c^\infty(\Omega \setminus \{0\})$ we set

$$J_n[v;\Omega] := \int_\Omega X^{-n+1}(|x|/R_\Omega)|\nabla v|^n dx.$$

The connection between $I_n$ and $J_n$ is demonstrated in the following proposition

**Proposition 2.6.** *Whenever $u \in C_c^\infty(\Omega \setminus \{0\})$ we have*

$$I_n[u;\Omega] \geq C_1(n) J_n[v;\Omega]; \quad C_1(n) = 1/(2^{n-1} - 1), \tag{2.4}$$

*where $v(x) := X^{1-1/n}(|x|/R_\Omega) u(x)$.*

**Proof.** Setting $u = X^{-1+1/n} v$ we compute

$$|\nabla u|^n = \left| X^{-1+1/n} \nabla v - \frac{n-1}{n} X^{1/n} \frac{v}{|x|} \frac{x}{|x|} \right|^n.$$

Now use the inequality $|b-a|^n \geq |a|^n + C_1(n)|b|^n - n|a|^{n-2} a \cdot b$ where $C_1(n) = 1/(2^{n-1} - 1)$ (see [7, Lemma 3.1]), to get

$$|\nabla u|^n \geq \left(\frac{n-1}{n}\right)^n \frac{|v|^n}{|x|^n} X + C_1(n) X^{-n+1} |\nabla v|^n - \left(\frac{n-1}{n}\right)^{n-1} \frac{1}{|x|^{n-1}} \nabla(|v|^n) \cdot \frac{x}{|x|}.$$

This means

$$\int_\Omega |\nabla u|^n dx \geq \left(\frac{n-1}{n}\right)^n \int_\Omega \frac{|v|^n}{|x|^n} X(|x|/R_\Omega) dx + C_1(n) J_n[v;\Omega]$$

$$+ \left(\frac{n-1}{n}\right)^{n-1} \int_\Omega |v|^n \operatorname{div}\left\{\frac{1}{|x|^{n-1}} \frac{x}{|x|}\right\} dx.$$

Since $\operatorname{div}\{\frac{1}{|x|^{n-1}} \frac{x}{|x|}\} = 0$ in $\Omega \setminus \{0\}$, we deduce

$$\int_\Omega |\nabla u|^n dx \geq \left(\frac{n-1}{n}\right)^n \int_\Omega \frac{|v|^n}{|x|^n} X(|x|/R_\Omega) dx + C_1(n) J_n[v;\Omega],$$

and rewriting the first term on the right with the original function $u$, we obtain (2.4). Note that for $n = 2$ we have equality in (2.4). ∎

## 3. Estimates in $L^q(\Omega)$; $q > n$

The main ingredient in the proof of (1.7) in Theorem 1.1 is the following estimate

**Proposition 3.1.** *For all $u \in C_c^\infty(\Omega \setminus \{0\})$, all $\varepsilon > 0$ and any $q > n \geq 2$, we have the inequality*

$$\|uX^\varepsilon\|_{L^q(\Omega)} \leq C_2(n,\varepsilon) \left(1 + q\frac{n-1}{n}\right)^{1/q + 1 - 1/n} \operatorname{vol}(\Omega)^{1/q} \left(J_n[v;\Omega]\right)^{1/n}, \tag{3.1}$$

*where $v(x) = u(x) X^{1-1/n}(|x|/R_\Omega)$, and $C_2(n,\varepsilon) = \frac{1}{n\omega_n^{1/n}}\left(1 + |\frac{n-1}{n\varepsilon} - 1|\right)$.*

**Proof.** It suffices to prove (3.1) for $R_\Omega = 1$. The general case follows by a change of variables. Setting $u(x) = v(x) X^{-1+1/n}(|x|)$, and then using the standard representation formula (see [18]-Lemma 7.14) we have

$$
\begin{aligned}
n\omega_n u(x) X^\varepsilon(|x|) &= n\omega_n v(x) X^{-1+1/n+\varepsilon}(|x|) \\
&= \int_\Omega \frac{(x-y) \cdot \nabla[v(y) X^{-1+1/n+\varepsilon}(|y|)]}{|x-y|^n} dy \\
&\leq \underbrace{\int_\Omega \frac{|\nabla v(y)| X^{-1+1/n}(|y|)}{|x-y|^{n-1}} dy}_{=:K(x)} + \left| \frac{n-1}{n} - \varepsilon \right| \underbrace{\int_\Omega \frac{|v(y)| X^{1/n+\varepsilon}(|y|)}{|y||x-y|^{n-1}} dy}_{=:\Lambda(x)},
\end{aligned}
$$

where in obtaining $K(x)$ we have used the fact that $X^\varepsilon(|x|) \leq 1$; $x \in \Omega$. Hence, by Minkowski's inequality

$$\|uX^\varepsilon\|_{L^q(\Omega)} \leq \frac{1}{n\omega_n}\left(\|K\|_{L^q(\Omega)} + \left|\frac{n-1}{n} - \varepsilon\right| \|\Lambda\|_{L^q(\Omega)}\right). \tag{3.2}$$

To bound $\|K\|_{L^q(\Omega)}$, we start by estimating K$(x)$. For any $q > n$ and $1 < r < n/(n-1)$ satisfying

$$\frac{1}{r} = \frac{n-1}{n} + \frac{1}{q}, \tag{3.3}$$

we may write the integrand of K$(x)$ as follows

$$\left\{\frac{1}{|x-y|^{(n-1)(1-r/q)}}\right\}\left\{|\nabla v(y)|^{1-n/q} X^{-(n-1)(1/n-1/q)}(|y|)\right\}\left\{\frac{|\nabla v(y)|^{n/q} X^{-(n-1)/q}(|y|)}{|x-y|^{(n-1)r/q}}\right\}.$$

Applying Hölder's inequality with conjugate exponents

$$\frac{n-1}{n} + \frac{q-n}{nq} + \frac{1}{q} = 1, \tag{3.4}$$

we obtain

$$K(x) \leq \|\mathscr{V}_r\|_{L^\infty(\Omega)}^{1-1/n} \left(J_n[v;\Omega]\right)^{1/n-1/q} \left(\int_\Omega \frac{|\nabla v(y)|^n X^{-n+1}(|y|)}{|x-y|^{(n-1)r}} dy\right)^{1/q},$$

where we have set $\mathscr{V}_r(x) := \int_\Omega |x-y|^{-(n-1)r} dy$. Taking the $L^q(\Omega)$ norm of the two sides we arrive at

$$
\begin{aligned}
\|K\|_{L^q(\Omega)} &\leq \|\mathscr{V}_r\|_{L^\infty(\Omega)}^{1-1/n} \left(J_n[v;\Omega]\right)^{1/n-1/q} \left(\int_\Omega \int_\Omega \frac{|\nabla v(y)|^n X^{-n+1}(|y|)}{|x-y|^{(n-1)r}} dy dx\right)^{1/q} \\
&= \|\mathscr{V}_r\|_{L^\infty(\Omega)}^{1-1/n} \left(J_n[v;\Omega]\right)^{1/n-1/q} \left(\int_\Omega |\nabla v|^n X^{-n+1}(|y|) \mathscr{V}_r(y) dy\right)^{1/q},
\end{aligned}
$$

by Tonelli's theorem. The last factor can be estimated by

$$\|\mathscr{V}_r\|_{L^\infty(\Omega)}^{1/q} \left(J_n[v;\Omega]\right)^{1/q},$$

and so

$$\|K\|_{L^q(\Omega)} \leq \|\mathscr{V}_r\|_{L^\infty(\Omega)}^{1/r} \left(J_n[v;\Omega]\right)^{1/n}. \tag{3.5}$$

Next we estimate $\Lambda(x)$ in order to obtain the analogous bound for $\|\Lambda\|_{L^q(\Omega)}$. The integrand of $\Lambda(x)$ can be written as follows

$$\left\{\frac{1}{|x-y|^{(n-1)(1-r/q)}}\right\}\left\{\frac{|v(y)|^{1-n/q}}{|y|^{1-n/q}}X^{(1+n\varepsilon)(1/n-1/q)}(|y|)\right\}\left\{\frac{|v(y)|^{n/q}X^{(1+n\varepsilon)/q}(|y|)}{|y|^{n/q}|x-y|^{(n-1)r/q}}\right\}.$$

Performing Hölder's inequality with the conjugate exponents (3.4), we get

$$\Lambda(x) \leq \|\mathscr{V}_r\|_{L^\infty(\Omega)}^{1-1/n}\left(\int_\Omega \frac{|v|^n}{|y|^n}X^{1+n\varepsilon}(|y|)dy\right)^{1/n-1/q}\left(\int_\Omega \frac{|v(y)|^n X^{1+n\varepsilon}(|y|)}{|y|^n|x-y|^{(n-1)r}}dy\right)^{1/q}.$$

Taking the $L^q(\Omega)$ norm of the two sides we arrive at

$$\|\Lambda\|_{L^q(\Omega)} \leq \|\mathscr{V}_r\|_{L^\infty(\Omega)}^{1-1/n}\left(\int_\Omega \frac{|v|^n}{|y|^n}X^{1+n\varepsilon}(|y|)dy\right)^{1/n-1/q}\left(\int_\Omega\int_\Omega \frac{|v(y)|^n X^{1+n\varepsilon}(|y|)}{|y|^n|x-y|^{(n-1)r}}dydx\right)^{1/q}$$

$$= \|\mathscr{V}_r\|_{L^\infty(\Omega)}^{1-1/n}\left(\int_\Omega \frac{|v|^n}{|y|^n}X^{1+n\varepsilon}(|y|)dy\right)^{1/n-1/q}\left(\int_\Omega \frac{|v|^n}{|y|^n}X^{1+n\varepsilon}(|y|)\mathscr{V}_r(y)dy\right)^{1/q},$$

by Tonelli's theorem. The last factor can be estimated by

$$\|\mathscr{V}_r\|_{L^\infty(\Omega)}^{1/q}\left(\int_\Omega \frac{|v|^n}{|y|^n}X^{1+n\varepsilon}(|y|)dy\right)^{1/q},$$

and so

$$\|\Lambda\|_{L^q(\Omega)} \leq \|\mathscr{V}_r\|_{L^\infty(\Omega)}^{1/r}\left(\int_\Omega \frac{|v|^n}{|y|^n}X^{1+n\varepsilon}(|y|)dy\right)^{1/n}$$

$$\leq \frac{1}{\varepsilon}\|\mathscr{V}_r\|_{L^\infty(\Omega)}^{1/r}\left(\int_\Omega |\nabla v|^n X^{-n+1+n\varepsilon}(|y|)dy\right)^{1/n}$$

$$\leq \frac{1}{\varepsilon}\|\mathscr{V}_r\|_{L^\infty(\Omega)}^{1/r}\left(J_n[v;\Omega]\right)^{1/n}, \tag{3.6}$$

where we have used first Lemma 2.1 with $\alpha = 1+n\varepsilon$ and then the fact that $X^{n\varepsilon}(|y|) \leq 1$ for all $y \in \Omega$. Inserting (3.6) and (3.5) in (3.2) we arrive at

$$\|uX^\varepsilon\|_{L^q(\Omega)} \leq \frac{1}{n\omega_n}\left(1+\left|\frac{n-1}{n\varepsilon}-1\right|\right)\|\mathscr{V}_r\|_{L^\infty(\Omega)}^{1/r}\left(J_n[v;\Omega]\right)^{1/n} \tag{3.7}$$

$$\leq \frac{1}{n\omega_n^{1/n}}\left(1+\left|\frac{n-1}{n\varepsilon}-1\right|\right)\omega_n^{1/q}\left(1+q\frac{n-1}{n}\right)^{1/q+1-1/n}\left(J_n[v;\Omega]\right)^{1/n},$$

by the fact that $\|\mathscr{V}_r\|_{L^\infty(\Omega)} \leq \frac{n\omega_n}{n-(n-1)r}$ and (3.3). ∎

## 4. Proof of Theorem 1.1

Let $u \in C_c^\infty(\Omega \setminus \{0\})$ be such that $I_n[u;\Omega] \leq 1$. By Proposition 2.6 we have $J_n[v;\Omega] \leq C_1^{-1}(n)$, where $v(x) = X^{1-1/n}(|x|/R_\Omega)u(x)$ and $C_1(n)$ is given in (2.4). Applying Proposition 3.1 with $q = ns/(n-1)$; $s \in \{n, n+1, ...\}$ we obtain

$$\frac{1}{\text{vol}(\Omega)} \int_\Omega \left[\left(|u(x)|X^\varepsilon(|x|)\right)^{n/(n-1)}\right]^s dx \leq \left(\frac{C_2(n,\varepsilon)}{C_1^{1/n}(n)}\right)^{sn/(n-1)} (s+1)^{s+1}.$$

Multiplying both sides by $c^s/s!$ and adding for all integers $s \in [n,k]$; $n \leq k \in \mathbb{N}$, gives

$$\int_\Omega \sum_{s=n}^k \frac{1}{s!}\left[c\left(|u(x)|X^\varepsilon(|x|)\right)^{n/(n-1)}\right]^s dx \leq \sum_{s=n}^k c^s \left(\frac{C_2(n,\varepsilon)}{C_1^{1/n}(n)}\right)^{sn/(n-1)} \frac{(s+1)^{s+1}}{s!}, \quad (4.1)$$

for any $k \in \{n, n+1, ...\}$, and $c > 0$ chosen so that the sum on the right hand converges as $k \to \infty$. It is enough to choose

$$c < \frac{1}{e}\left(\frac{C_1^{1/n}(n)}{C_2(n,\varepsilon)}\right)^{n/(n-1)}.$$

Using Jensen's inequality and then Proposition 3.1 we see that each term of the finite sum

$$S = \frac{1}{\text{vol}(\Omega)} \int_\Omega \sum_{s=0}^{n-1} \frac{1}{s!}\left[c\left(|u(x)|X^\varepsilon(|x|)\right)^{n/(n-1)}\right]^s dx,$$

is bounded by a constant that depends only on $n$. Thus, adding $S$ on both sides of (4.1), the proof in case of $\Omega$ with $R_\Omega = 1$ is completed by letting $k \to \infty$ and using the monotone convergence theorem. The case of general $\Omega$ follows by scaling.

Next we show that (1.7) fails for $\varepsilon = 0$. Suppose, for the sake of contradiction, that there exist positive constants $c_1, c_2$ such that

$$\int_{B_1} e^{c_1|u(x)|^{n/(n-1)}} dx \leq c_2 \quad \text{for all } u \in C_c^\infty(B_1 \setminus \{0\}) \text{ satisfying } I_n[u;B_1] \leq 1. \quad (4.2)$$

Then we claim that one can obtain the inequality

$$\int_{B_1} \frac{|u|^n}{|x|^n} X^n(|x|) X^\theta(X(|x|)) dx \leq C I_n[u;B_1] \quad \text{for all } u \in C_c^\infty(B_1 \setminus \{0\}), \quad (4.3)$$

for any $\theta > 1$, an absurdity, since this is only possible if $\theta \geq 2$ (see [8, Theorem B]).

We therefore proceed to establish inequality (4.3). Let us recall the following version of Young's inequality

$$ts \leq e^t - t - 1 + (1+s)\log(1+s) - s, \quad (4.4)$$

for $s, t \geq 0$. Here we have taken the conjugate functions $A(t) = e^t - t - 1$ and $\tilde{A}(s) = (1+s)\log(1+s) - s$ (see for instance [1, §8.3]). Now, we set $t = a^{1/(n-1)}$ and $s = b^{1/(n-1)}$ in (4.4) and arrive to the inequality

$$\begin{aligned} ab &\leq \left(e^{a^{1/(n-1)}} + (1+b^{1/(n-1)})\log(1+b^{1/(n-1)})\right)^{n-1} \\ &\leq 2^{n-2}\left(e^{(n-1)a^{1/(n-1)}} + 2^{n-2}(1+b)[\log(1+b^{1/(n-1)})]^{n-1}\right) \quad \text{for all } a,b \geq 0. \end{aligned} \quad (4.5)$$

Writing

$$\int_{B_1} \frac{|u|^n}{|x|^n} X^n(|x|) X^\theta(X(|x|)) dx = \left(\frac{n-1}{c_1}\right)^{n-1} \int_{B_1} \left\{\left(\frac{c_1}{n-1}\right)^{n-1} |u|^n\right\} \left\{\frac{X^n(|x|) X^\theta(X(|x|))}{|x|^n}\right\} dx,$$

and applying (4.5), we have, because of (4.2), that

$$\int_{B_1} \frac{|u|^n}{|x|^n} X^n(|x|) X^\theta(X(|x|)) dx \leq 2^{n-2}\left(\frac{n-1}{c_1}\right)^{n-1} c_2 + M, \quad (4.6)$$

where

$$M := 2^{2(n-2)} \int_{B_1} \left(1 + \frac{X^n(|x|) X^\theta(X(|x|))}{|x|^n}\right) \left[\log\left(1 + \left(\frac{X^n(|x|) X^\theta(X(|x|))}{|x|^n}\right)^{1/(n-1)}\right)\right]^{n-1} dx.$$

It remains to demonstrate that $M$ is finite for $\theta > 1$, since (4.6) will then imply

$$\int_{B_1} \frac{|u|^n}{|x|^n} X^n(|x|) X^\theta(X(|x|)) dx \leq C \quad \text{for all } u \in C_c^\infty(B_1 \setminus \{0\}) \text{ satisfying } I_n[u; B_1] \leq 1,$$

which by the normalization

$$\tilde{u} = \frac{u}{\left(I_n[u, B_1]\right)^{1/n}},$$

yields (4.3).

Now, since the integrand in $M$ is bounded away from 0, it suffices to estimate the integral in a small ball around the origin. Let $\eta > 0$ be such that $|x|^{-n} X^n(|x|) X^\theta(X(|x|)) \geq 1$ for any $x \in B_\eta$, and let $M_\eta$ denote the integral on $B_\eta$. We estimate

$$M_\eta \leq 2^{2n-3} \int_{B_\eta} \frac{X^n(|x|) X^\theta(X(|x|))}{|x|^n} \left[\log\left(2\left(\frac{X^n(|x|) X^\theta(X(|x|))}{|x|^n}\right)^{1/(n-1)}\right)\right]^{n-1} dx,$$

and using polar coordinates

$$M_\eta \leq n\omega_n 2^{2n-3} \int_0^\eta t^{-1} X^n(t) X^\theta(X(t)) \left[\log\left(2t^{-n/(n-1)} X^{n/(n-1)}(t) X^{\theta/(n-1)}(X(t))\right)\right]^{n-1} dt.$$

Now notice that $\log\left(2t^{-n/(n-1)} X^{n/(n-1)}(t) X^{\theta/(n-1)}(X(t))\right) \leq \frac{n}{n-1} X^{-1}(t)$ to conclude

$$M_\eta \leq n\omega_n 2^{2n-3} \left(\frac{n}{n-1}\right)^{n-1} \int_0^\eta t^{-1} X(t) X^\theta(X(t)) dt.$$

This last integral is finite if and only if $\theta > 1$ (see for example [8, Proposition 3.1-Equation (3.8)]), which proves the claim and yields the desired contradiction. ∎

**Acknowledgements** The first author is supported in part at the Technion by a Fine Fellowship and the second author is supported in part by a Technion fellowship.